\documentclass{article}

\usepackage{graphicx}
\graphicspath{ {./images/} }

\begin{document}
\title{\centering \large \bf
New Isothermic surfaces }
\author{\large \bf
Armando M. V. Corro\\
IME, Universidade Federal de Goi\'as \\
Caixa Postal 131, 74001-970, Goi\^ania, GO, Brazil\\
e-mail:corro@ufg.br
  \and
 \large \bf Marcelo Lopes Ferro    \\
 IME, Universidade Federal de Goi\'as \\
Caixa Postal 131, 74001-970, Goi\^ania, GO, Brazil\\
e-mail: marceloferro@ufg.br
 }
\date{}
\maketitle \thispagestyle{empty}
 {\begin{abstract}
In this paper, we consider a method of constructing isothermic
surfaces based on Ribaucour transformations. By applying the theory to the cylinder, we obtain a three-parameter family of complete isothermic surfaces that contains n-bubble surfaces inside and outside of the cylinder. In addition, we also obtain one-parameter family of complete isothermic surface with planar ends. Such family of isothermic surfaces do not have constant mean curvature. As aplication we obtain explicit solutions of the Calapso equation. 
\end{abstract}

\section*{\large \bf  Introduction}

Ribaucour transformations for hypersurfaces, parametrized by lines of curvature, were classically studied by Bianchi \cite{1}. They can be applied to obtain surfaces of constant Gaussian curvature and surfaces of constant mean curvature,  from a given such surface, respectively, with constant Gaussian curvature and constant mean curvature. The first application of this method to minimal and
cmc surfaces in $R^3$ was obtained by Corro, Ferreira, and Tenenblat in \cite{7}-\cite{8}. For more application this method, see \cite{9}-\cite{10}, \cite{17}, \cite{22}, \cite{23} and \cite{24}.

Using Ribaucour transformations and applying the theory to the cylinder, \cite{8} obtained a two-parameter family of complete linear Weingarten surfaces, that contains n-bubble surfaces inside of the cylinder.

A regular surface $M$ is isothermic if locally, near each non umbilic point of M there exist curvature line coordinates which are conformal with respect to the first fundamental form of $M$.

The study of isothermic surfaces is a very difficult problem because it depends on the integration of an equation with partial derivatives of fourth order ( see \cite{25} ). Particular classes of these surfaces can be obtained using some transformations from a given isothermic surface.

The theory of isothermic surfaces was studied by eminent  geometries as Christoffel \cite{5}, Darboux \cite{15}-\cite{16} and Bianchi \cite{1} among others. Particular classes of isothermic surfaces are, the constant mean curvature surfaces, quadrics surfaces, surfaces whose lines of curvature has constant geodesic curvature, in particular, the cyclides of Dupin.  The isothermic surfaces are preserved by isometries, dilations and invertions.

In \cite{2} the authors study surfaces with harmonic inverse mean curvature (HIMC surfaces). They distinguish a subclass of $\theta-isothermic$ susrfaces, and if $\theta=0$ then surfaces are isothermic.

In \cite{4} the author show that theory of soliton surfaces, modified in an approprieate way, can be applied also to isothermic immersions in $R^3$. In this case the so called Sym's formula gives an explicit expression for the isothermic immersion with prescribed fundamental forms. The complete classification of the isothermic surfaces is an open problem.

In \cite{3}, the author establishes an equation with fourth order partial derivatives from which the problem of obtaining isothermic surfaces apparently becomes much simpler. Such equation ( called Calapso equation ) defined in \cite{3}  given by
\begin{eqnarray}\nonumber
\bigg(\frac{\omega,_{u_1u_2}}{w}\bigg),_{u_1u_1}+\bigg(\frac{\omega,_{u_1u_2}}{w}\bigg),_{u_2u_2}+\big(\omega^2\big),_{u_1u_2}=0,\nonumber
\end{eqnarray}
describes isothermic surfaces in $R^3$.

In \cite{13} the outhors introduced the class of radial inverse mean curvature survace ( RIMC - surfaces ), that are isothermics surfaces. Moreover, were obtained two solutions of the Calapso equation where one can be obtained using \cite{3} and a different one.

In this paper, motivated by \cite{8}  and \cite{13} we use the  Ribaucour transformations to get a family of isothermic surfaces from a given such surface. As an application of the theory, we obtain a family of complete isothermic surfaces associated to the cylinder  with no constant mean curvature. 

The families we obtain depend on one or three parameters. One of them is the parameter $c\neq 0$ of the Ribaucour transformation. According to the value of c, we can have n-bubbles surface or not. Precisely, if $\sqrt{|1+c|}=n/m$ an irreducible rational number, we get an n-bubble surface. This is an immersed cylinder into $R^3$, with two ends of geometric index $m$ and $n$ isolated points of maximum and of minimum for the Gaussian curvature. Otherwise, if $\sqrt{|1+c|}$ is not a rational number,  then isothermic surface is a complete immersion of $R^2$ into $R^3$, not periodic in any variable. Another parameter is $b$, and depending on the b sign, we have n-bubbles surface inside or outside of the cylinder. The last parameter, appear from integrating the Ribaucour transformation. In the familie depending on one-parameter,  we have planar ends. The position of the planar ends is affect by the value of the parameter. We show explicit examples that check for planar ends and bubbles surfaces.

Also, in this work motivated by \cite{13}, for each isothermic surfaces obtained by the Ribaucour transformations, we associates a solution of the Calapso equation. We give explicit solutions of the Calapso equation that depend on functions of a single variable. Applying isometries, dilations,
inversions,  we  obtain new isothermic surfaces. So we get new solutions of the Calapso equation.

 \vspace{.1in}

\section*{ \large \bf 1. Ribaucour transformation for isothermic surfaces}
\vspace{.1in}

In this section, we first recall the theory of Ribaucour
transformation for surfaces ( see \cite{1} and \cite{7}  for more details ).

A isothermic surface of $R^3$ is a surface which locally has an
orthogonal coordinate system with same coefficients of the first
quadratic form.

Let $M$ be an orientable surface of $R^3$ without umbilic points,
with Gauss map we denote by $N$. We say that $\widetilde{M}$ is
associated to $M$ by a Ribaucour transformation, if and only if,
there exists a differentiable function $h$ defined on $M$ and a
diffeomorphism $\psi:M\rightarrow \widetilde{M}$ such that\\
(a) for all $p\in M$,
$p+h(p)N(p)=\psi(p)+h(p)\widetilde{N}(\psi(p))$, where
$\widetilde{N}$ is the Gauss map of $\widetilde{M}$.\\
(b)The subset $p+h(p)N(p)$, $p\in M$, is a two-dimensional
submanifold.\\
(c) $\psi$ preserves lines of curvature.

We say that $\widetilde{M}$ is locally associated to $M$ by a
\textit{Ribaucour transformation} if, for all $\widetilde{p}$,
there exists a neighborhood of $\widetilde{p}$ in $\widetilde{M}$
which is associated by a Ribaucour transformation to an open
subset of $M$.

The following result gives a characterization of Ribaucour
transfomations. For the proof and more details, see \cite{8}

\vspace{.1in}

\noindent {\bf Theorem 1.1.} \textit{Let $M$ be an orientable
surface of $R^3$, without umbilic points, whose Gauss map is $N$.
Let $e_i$, $1\leq i\leq 2$ be orthonormal principal directions,
$-\lambda_i$ the corresponding principal curvatures, i.e.
$dN(e_i)=\lambda_ie_i$. A surface $\widetilde{M}$ is associated to
$M$ by a Bibaucour transformation, if and only if, $M$ and
$\widetilde{M}$ are associated by a sphere congruence whose radius
function $h:M\rightarrow R$ satisfies $1+h\lambda_i\neq 0$ and
\begin{equation}\label{eq1}
dZ^j(e_i)+Z^iw_{ij}-Z^iZ^j\lambda_i=0,\hspace{1cm} 1\leq i\neq
j\leq 2,
\end{equation}
where
\begin{eqnarray*}\nonumber
Z^i=\frac{dh(e_i)}{1+h\lambda_i},\nonumber
\end{eqnarray*}
and $w_{ij}$ are the connection forms of the frame $e_i$.
 }

 \vspace{.1in}

 \noindent {\bf Remark 1.2} In a local coordinate system by lines
 of curvature, the function $h(u_1,u_2)$ is a differentiable
 function which satisfies a second-order nonlinear partial
 differential equation corresponding to (\ref{eq1}). One can
 linearize the problem of obtaining the function $h$. This is a
 consequence of the following result ( see \cite{9}
for a proof and more details)

 \vspace{.1in}

 \noindent {\bf Proposition 1.2} \textit{If $h$ is a nonvanishing
 function which satisfies (\ref{eq1}), then
\begin{eqnarray*}\nonumber
\frac{1}{h}\sum_{i=1}^2Z^iw_i,\nonumber
\end{eqnarray*}
is a closed $1-form$ and there is a nonvanishing function
$\Omega$, defined on a simply conected domain, such that
\begin{eqnarray*}\nonumber
d\Omega(e_i)=\frac{\Omega}{h}Z^i.\nonumber
\end{eqnarray*}}

 \vspace{.1in}

For each nonvanishing function $h$, which is a solution of
(\ref{eq1}), we consider $\Omega$ as above and we define
\begin{eqnarray*}\nonumber
\Omega_i=d\Omega(e_i),\hspace{1cm}W=\frac{\Omega}{h}.\nonumber
\end{eqnarray*}
With this notation,
\begin{eqnarray}
dh(e_i)=\frac{\Omega_i}{W}\bigg(1+\frac{\Omega\lambda_i}{W}\bigg),\hspace{0,5cm}1+h\lambda_i=1+\frac{\Omega\lambda_i}{W},\hspace{0,5cm}Z^i=\frac{\Omega_i}{W}.\label{eq2}
\end{eqnarray}

\vspace{.1in}

The next result shows that (\ref{eq1}) is equivalent to a linear
system, for more details, see \cite{8}.

 \vspace{.1in}

\noindent {\bf Proposition 1.4} \textit{ A function $h$ is a
solution of (\ref{eq1}) defined on a simply connected domain, if
and only if, $h=\frac{\Omega}{W}$, where $\Omega$ and $W$ are
functions which satisfy
\begin{eqnarray}
d\Omega_i(e_i)&=&\Omega_jw_{ij}(e_j),\hspace{0,5cm}
for\hspace{0,2cm}i\neq j, \label{eq3}\\
d\Omega &=&\sum_{i=1}^n\Omega_iw_i,\label{eq4}\\
dW&=&-\sum_{i=1}^n\Omega_i\lambda_iw_i.\label{eq5}
\end{eqnarray}
}

 \vspace{.1in}

One can show that the Ribaucour transformation of a surface is
given in terms of the solution of the above system ( see \cite{7}
for a proof and more details)

\vspace{.1in}

\noindent {\bf Theorem 1.5} \textit{ Let $M$ be an orientable
surface of $R^3$ parametrized by $X:U\subseteq R^2\rightarrow M$,
without umbilic points. Assume $e_i$, $1\leq i\leq 2$ are
orthogonal principal directions, $-\lambda_i$ the corresponding
principal curvatures, and $N$ is a unit vector field normal to
$M$. A surface $\widetilde{M}$ is locally associated to $M$ by a
Ribaucour transformation if and only if there is differentiable
functions $W,\Omega,\Omega_i:V\subseteq U\rightarrow R$ which
satisfy (\ref{eq3}), (\ref{eq4}), (\ref{eq5}),
$W(W+\lambda_i\Omega)\neq 0$ and $\widetilde{X}:V\subseteq
U\rightarrow \widetilde{M}$, is a parametrization of
$\widetilde{M}$ given by
\begin{eqnarray}
\widetilde{X}=X-\frac{2\Omega}{S}\bigg(\sum_{i=1}^2\Omega_ie_i-WN\bigg),\label{eq6}
\end{eqnarray}
where
\begin{eqnarray}
S=\sum_{i=1}^2\big(\Omega_i\big)^2+W^2.\label{eq7}
\end{eqnarray}
Moreover, the normal map of $\widetilde{X}$ is given by
\begin{eqnarray}
\widetilde{N}=N+\frac{2W}{S}\bigg(\sum_{i=1}^2\Omega_ie_i-WN\bigg),\label{eq8}
\end{eqnarray}
and the principal curvatures and coefficients of the first
fundamental form of $\widetilde{X}$, are given by
\begin{eqnarray}
&&\widetilde{\lambda}_i=\frac{WT_i+\lambda_iS}{S-\Omega
T_i},\label{eq9}\\
&&\widetilde{g}_{ii}=\bigg(\frac{S-\Omega
T_i}{S}\bigg)^2g_{ii}\label{eq10}
\end{eqnarray}
where $\Omega_i$, $\Omega$ and $W$ satisfy (\ref{eq3}),
(\ref{eq4}), (\ref{eq5}), $S$ is given by (\ref{eq7}), $g_{ii}$,
$1\leq i\leq 2$ are coefficients of the first fundamental form of
$X$, and
\begin{eqnarray}
T_i=2\bigg(d\Omega_i(e_i)+\sum_{k}\Omega_kw_{ki}(e_i)-W\lambda_i\bigg),\label{eq11}
\end{eqnarray}
}

\vspace{.1in}

\noindent {\bf Remark 1.6} If M is parametrized by orthogonal
curvature lines, we will assume that $e_i$ are given by
$\frac{X,_i}{a_i}$, where $a_i=\sqrt{g_{ii}}$. In this case,
system (\ref{eq3}), (\ref{eq4}) and (\ref{eq5}) can be rewritten
by
\begin{eqnarray}
\Omega_{i,j}&=&\Omega_j\frac{a_j}{a_{j,i}},\hspace{0,5cm}
for\hspace{0,2cm}i\neq j, \nonumber\\
\Omega,_i &=&a_i\Omega_i,\label{eq12}\\
W,_i&=&-a_i\Omega_i\lambda_i.\nonumber
\end{eqnarray}

\vspace{.1in}

The two result following provides a sufficient condition for a
Ribaucour transformation to transform a isothermic surface into
another such surface. The first one is

\vspace{.1in}

\noindent {\bf Theorem 1.7} \textit{Let $M$ be a surfaces of $R^3$
parametrized by $X:U\subseteq R^2\rightarrow M$, without umbilic
points and let $\widetilde{M}$  parametrized by (\ref{eq6}) be
associated to $M$ by a Ribaucour transformation, such that the
normal lines intersect at a distance function $h$. Assume that
$h=\frac{\Omega}{W}$ is not constant along the lines of curvature
and the function $\Omega$, $\Omega_i$ and $W$ satisfy one of the
additional relation
\begin{eqnarray}
 T_1+T_2=\frac{2S}{\Omega}\hspace{0,5cm}or\hspace{0,5cm}T_1-T_2=0\label{eq13}
\end{eqnarray}
where $S$ is given by (\ref{eq7}) and $T_i$, $1\leq i\leq 2$ are
defined by (\ref{eq11}). Then $\widetilde{M}$ parameterized by
(\ref{eq6}) is a isothermic surface, if and only if $M$ is
isothermic surface.
 }

\vspace{.1in}

\noindent {\bf Proof:}Suppose that $\widetilde{M}$ is a isothermic
surface, then the coefficients of the first fundamental form of
$\widetilde{X}$ satisfy, $\widetilde{g}_{11}=\widetilde{g}_{22}$.
So, using (\ref{eq10}), we have
\begin{eqnarray}
 \bigg(\frac{S-\Omega
T_1}{S}\bigg)^2g_{11}=\bigg(\frac{S-\Omega
T_2}{S}\bigg)^2g_{22},\label{eq14}
\end{eqnarray}
where $g_{ii}$, $1\leq i\leq 2$ are the coefficients of the first
fundamental form of $X$.\\
If $T_1+T_2=\frac{2S}{\Omega}$, then isolating $T_1$ and
substituting in (\ref{eq14}), we get $g_{11}=g_{22}$.\\
On the other hand, if $T_1-T_2=0$, then we have from (\ref{eq14})
that $g_{11}=g_{22}$. Therefore, $M$ is a isothermic surface.

Conversely, suppose that $M$ is a isothermic surface, then using
(\ref{eq13}), immediately from (\ref{eq10}), we obtain that
$\widetilde{M}$ is a isothermic surface.

\vspace{.1in}

 \noindent {\bf Remark 1.8}  Let $X:U\subseteq R^2\rightarrow M$ a
 isothermic parametrized for $M$. So, the first fundamental form of
 $X$ is given by $I=e^{2\varphi}\big(du_1^2+du_2^2\big)$. Thus,
 the first additional relation of (\ref{eq13}) is equivalent to
\begin{eqnarray}
\Delta\Omega-e^{2\varphi}(\lambda_1+\lambda_2)W=\frac{Se^{2\varphi}}{\Omega}\label{eq15}
\end{eqnarray}
In fact, under these conditions $T_i$ $1\leq i\leq 2$, given
 by (\ref{eq11}), can be rewritten as
\begin{eqnarray*}
&&T_1=\frac{2}{e^{\varphi}}\bigg(\Omega_{1,1}+\varphi,_{2}\Omega,_2-W\lambda_1e^{\varphi}\bigg),\nonumber\\
&&T_2=\frac{2}{e^{\varphi}}\bigg(\Omega_{2,2}+\varphi,_{1}\Omega,_1-W\lambda_2e^{\varphi}\bigg)\nonumber
\end{eqnarray*}
Using (\ref{eq12}) in this last equation, we get
\begin{eqnarray*}
&&T_1=\frac{2}{e^{2\varphi}}\bigg(\Omega_{,11}-\varphi,_1\Omega,_1+\varphi,_{2}\Omega,_2-W\lambda_1e^{2\varphi}\bigg),\nonumber\\
&&T_2=\frac{2}{e^{2\varphi}}\bigg(\Omega_{,22}-\varphi,_2\Omega,_2+\varphi,_{1}\Omega,_1-W\lambda_2e^{2\varphi}\bigg).\nonumber
\end{eqnarray*}
Therefore,
$T_1+T_2=\frac{2}{e^{2\varphi}}\big(\Delta\Omega-e^{2\varphi}W(\lambda_1+\lambda_2)\big)$
and the first additional relation of (\ref{eq13}) is equivalent to
(\ref{eq15}).

\vspace{.2in}

 \noindent {\bf Remark 1.9}  Let $X$ as in the previous remark.
Then the parameterization $\widetilde{X}$ of $\widetilde{M}$,
locally associated to $X$ by a Ribaucour transformation, given by
(\ref{eq6}), is defined on
\begin{eqnarray*}
V=\{(u_1,u_2)\in U;\hspace{0,1cm}\Omega T_1-S\neq0\}.\nonumber
\end{eqnarray*}

\vspace{.2in}

\section*{ \large \bf 2. Families of isothermic surfaces associated to the cylinder.}

\vspace{.2in}

In this section, by applying Theorem 1.7, using the Remark 1.8 and
1.9 to the cylinder, we obtain a three-parameter family of
complete isothermic surfaces. As obtained in \cite{8}, using examples, we check that there are n-bubble surfaces outside which are 1-periodic,
have genus zero and two ends of finite geometric index. In addition, we  have 1-periodic n-bubble surfaces inside of the cylinder and  one-parameter family of complete isothermic surfaces with planar ends.

 \vspace{.2in}

 \noindent {\bf Theorem 2.1} \textit{ Consider the cylinder
 parametrized by
\begin{eqnarray}
X(u_1,u_2)=(\cos(u_2),\sin(u_2),u_1), \hspace{0,3cm}(u_1,u_2)\in
R^2\label{cilindro}
\end{eqnarray}
as isothermic surface where the first fundamental form is
$I=du_1^2+du_2^2$. A parametrized surface $\widetilde{X}(u_1,u_2)$
is isothermic surface locally associated to $X$ by a Ribaucour
transformation as in Theorem 1.7, if and only if, up to a rigid
motion of $R^3$, it is given by
\begin{eqnarray}
\widetilde{X}=X-\frac{2}{2b+c\big(f-g\big)}\bigg(f'X,_1+g'X,_2-gN\bigg)\label{cilindroribaucour}
\end{eqnarray}
defined on $V=\{(u_1,u_2)\in R^2;\hspace{0,1cm}f+g\neq0\}$ where
$N$ is the inner unit normal vector field of the cylinder,
$c\neq0$, $b$ is real constant, and $f(u_1)$, $g(u_2)$ are
solutions of the equations
\begin{eqnarray}
&&f''-cf=b,\label{eqf}\\
&&g''+(1+c)g=b\nonumber
\end{eqnarray}
with initial conditions satisfying
\begin{eqnarray}
\big[\big(f'\big)^2-cf^2-2bf+\big(g'\big)^2+(1+c)g^2-2bg\big](u_1^0,u_2^0)=0.\label{algcond}
\end{eqnarray}
Moreover, the normal map of $\widetilde{X}$ is given by
\begin{eqnarray}
\widetilde{N}=N-\frac{2g}{(2b+c\big(f-g\big))(f+g)}\bigg(f'X,_1+g'X,_2-gN\bigg)\label{normalcilindroribaucour}
\end{eqnarray}
}

\vspace{.1in}

\noindent {\bf Proof:} Consider the first fundamental form of the
cylinder $ds^2=du_1^2+du_2^2$ and the principal curvatures
$\lambda_1=0$, $\lambda_2=-1$. Using (\ref{eq12}), to obtain the
Ribaucour transformations, we need to solve the following of
equations
\begin{eqnarray}
\Omega_{i,j}=0,\hspace{0,5cm}\Omega,_i
=\Omega_i,\hspace{0,5cm}W,_i=-\Omega_i\lambda_i,
\hspace{0,2cm}1\leq i\neq j\leq 2.\nonumber
\end{eqnarray}

Since $\Omega,_{12}=0$, it follows that $\Omega =
f_1(u_1)+g_2(u_2)$. Therefore $\Omega_1=f_1'$ and $\Omega_2=g_2'$.
Moreover, $W=g_2+a$, where $a$ is a real constant. Thus, from
(\ref{eq7}),
$S=(f_1')^2+(g_2')^2+(g_2+a)^2$. \\
Using the Remark 1.8, the associated surface will be isothermic
when $\Delta\Omega + W=\frac{S}{\Omega}$. Therefore, we obtain
that the functions $f_1$ and $g_2$ satisfy
\begin{eqnarray}
f_1''+g_2''+g_2+a=\frac{(f_1')^2+(g_2')^2+(g_2+a)^2}{f_1+g_2}.\label{eqf1''}
\end{eqnarray}
Differentiate this last equation with respect $x_1$ and $x_2$, we
get
\begin{eqnarray}
&&f_1'''=f_1'\bigg(\frac{f_1''-g_2''-g_2-a}{f_1+g_2}\bigg).\label{eqf'''}\\
&&g_2'''+g_2'=-g_2'\bigg(\frac{f_1''-g_2''-g_2-a}{f_1+g_2}\bigg).\nonumber
\end{eqnarray}

 Differentiate $\displaystyle{\frac{f_1''-g_2''-g_2-a}{f_1+g_2}}$,
with respect $x_i$, $1\leq i\leq 2$ and using (\ref{eqf'''}), we
get $\displaystyle{\bigg(\frac{f_1''-g_2''-g_2-a}{f_1+g_2}\bigg)_i=0}$.
Therefore  $\displaystyle{\frac{f_1''-g_2''-g_2-a}{f_1+g_2}=c}$,
where $c$ is a real constant. Thus, we have that $f_1$ and $g_2$
satisfy

\begin{eqnarray}\nonumber
&&f_1''-cf_1+ca=b,\nonumber \\
&&g_2''+(1+c)g_2+a+ac=b.\nonumber
\end{eqnarray}
Now defining $f=f_1-a$ and $g=g_2+a$, we obtain that $f$ and $g$
satisfy (\ref{eqf}), with real constant $c\neq0$, because if
$c=0$, then $\widetilde{X}$ is degenerate. Moreover, using
(\ref{eqf1''}) we get that the initial conditions satisfying
(\ref{algcond}) and using the Theorem 1.5, $\widetilde{X}$ is
give by (\ref{cilindroribaucour}) and from Remark 1.9 is defined
in $V=\{(u_1,u_2)\in R^2;\hspace{0,1cm}f+g\neq 0\}$.


 \vspace{.2in}



 \noindent {\bf Remark 2.2} Each isothermic surfaces associated to
 the cylinder as in Theorem 2.1, is parametrized by lines of
 curvature and from (\ref{eq10}), the metric is given by
 $ds^2=\psi^2(du_1^2+du_2^2)$, where
\begin{eqnarray}
\psi=\frac{|c\big(f+g\big)|}{|2b+c\big(f-g\big)|}.\label{metrica}
\end{eqnarray}
Moreover, from (\ref{eq9}), the principal curvatures of the
$\widetilde{X}$ are given by
\begin{eqnarray}
\widetilde{\lambda}_1&=&\frac{-2g\big(b+cf\big)}{c\big(f+g\big)^2},\label{aa1}\\
\widetilde{\lambda}_2&=&\frac{-cf^2-2bf-cg^2}{c\big(f+g\big)^2}.\label{aa2}
\end{eqnarray}

\vspace{.2in}

\noindent {\bf Proposition 2.3} \textit{Consider the isothermic
surfaces associated to the cylinder parametrized by
(\ref{cilindroribaucour}). Then the mean curvature of the
$\widetilde{X}$ is given by
\begin{eqnarray}
\widetilde{H}=\frac{-1}{2}-\frac{b}{c\big(f+g\big)}.\label{mediaribo}
\end{eqnarray}
}

\vspace{.1in}

\noindent {\bf Proof:} In fact, using Remark 2.2 is easy to get
\begin{eqnarray}
\widetilde{\lambda}_1+\widetilde{\lambda}_2=\frac{-\big(cf+cg+2b\big)}{c\big(f+g\big)}.\nonumber
\end{eqnarray}
Therefore, $\widetilde{H}$ is given by (\ref{mediaribo}).


\vspace{.2in}

Using the previous proposition, we immediately get

\vspace{.1in}

\noindent {\bf Corollary 2.4} \textit{Consider the isothermic
surfaces associated to the cylinder parametrized by
(\ref{cilindroribaucour}). Then $\widetilde{X}$ is
$\frac{-1}{2}-cmc$, if and only if, $b=0$.}

\vspace{.2in}

The next result, describes the behavior of  $\widetilde{X}$, in the neighborhood of $p_0$ where\\
$\big(2b+c(f-g)\big)(p_0)=0$ .

\vspace{.2in}

\noindent {\bf  Proposition 2.5} \textit{Let $M=2b+c\big(f-g\big)$ and consider $p_0\in R^2$ such that $M(p_0)=0$ and $M(p)\neq 0$, for all $p\in V-\{p_0\}$. Let  $\widetilde{X}:V-\{p_0\}\rightarrow R^3$ be a isothermic surface locally associated by a Ribaucour 
transformation to  cylinder given by  Theorem 2.1, where $V=\{(u_1,u_2)\in R^2;\hspace{0,1cm}f+g\neq 0\}$. Then for any divergent curve $\gamma:[0,1)\rightarrow V-\{p_0\} $ such that $\lim_{t\rightarrow 1}\gamma(t)=p_0$ the length of $\widetilde{X}(\gamma)$
is infinite.}

\vspace{.2in}

\noindent {\bf Proof:} Let $p_0\in R^2$ such that $M(p_0)=0$, then using
(\ref{eqf}) and (\ref{eqf1''}), we have $S(p_0)=0$. After a
translation, we may assume that $p_0=(0,0)$. It follows from Remark 2.3 that the first
fundamental of $\widetilde{X}$ is given by
$\widetilde{I}=\psi^2(du_1^2+du_2^2)$, where
$\psi=\displaystyle{\frac{|c\big(f+g\big)|}{|M|}}$.

At $p_0$ we have $M(0,0)=0$ and $S(0,0)=0$, therefore
\begin{eqnarray}
f'(0)=0,\hspace{0,5cm}g'(0)=0,\hspace{0,5cm}g(0)=0 \hspace{0,5cm}and\hspace{0,5cm}f(0)=\displaystyle{\frac{-2b}{c}}.\label{fgponto}
\end{eqnarray}

Differentiate $M$ and evaluating and $(0,0)$, we get
$M,_i(0,0)=M,_{ij}(0,0)=0$ for $1\leq i\neq j\leq 2$ and
$M,_{ii}=-cb$. Hence, by considering the Taylor expansion of $M$ on
neighborhood of $(0,0)$, we obtain
\begin{eqnarray}
M(p)=\frac{-cb|p|^2}{2}+R,\hspace{0,3cm}with\hspace{0,3cm}\lim_{|p|\rightarrow
0}\frac{R}{|p|^2}=0,\hspace{0,2cm}p=(u_1,u_2).\label{Mu}
\end{eqnarray}
Therefore, we have
\begin{eqnarray}
\lim_{|p|\rightarrow 0}|p|^2\widetilde{\psi}=\frac{4}{c},\nonumber
\end{eqnarray}
where $\widetilde{\psi}=\displaystyle{\frac{c\big(f+g\big)}{M}}$.

Hence, there exists $\delta >0$ such that, $0<|p|<\delta$ we have
$|p|^2\psi>\displaystyle{\frac{2}{|c|}}$.\\
Let divergent curve $\gamma:[0,1)\rightarrow V-\{(0,0)\}$, where
$V=\{(u_1,u_2)\in R^2;\hspace{0,1cm}f+g\neq 0\}$, such that
$\displaystyle{\lim_{t\rightarrow 1}\gamma(t)=(0,0)}$. Then its
length is
\begin{eqnarray}
l(\widetilde{X}\circ\gamma)=\int_0^1\psi(\gamma(t))|\gamma'(t)|dt>\frac{2}{|c|}\int_0^1\frac{|\gamma'(t)|}{|\gamma(t)|^2}dt\geq
\bigg|\frac{2}{|c|}\int_0^1\frac{|\gamma'(t)|}{|\gamma(t)|^2}dt\bigg|=\infty.
\end{eqnarray}

\vspace{.2in}

As in \cite{8}, we introduce a notation, which will be useful in the
following results.\\
A \textit{rotation of angle} $\theta$ in the $xy$ plane of $R^3$
will be denoted by
\begin{equation}
R_\theta=\left[
\begin{array}{ccc}
  \cos \theta & -\sin \theta & 0 \\
  \sin \theta & \cos \theta & 0 \label{rotacao}\\
  0           &     0        & 1
\end{array}
\right].
\end{equation}
We denote by $T_\delta$ the translation defined by
\begin{eqnarray}\label{translacao}
T_\delta(x,y,z)=(x,y,z+\delta).
\end{eqnarray}

\vspace{.2in}

\noindent {\bf Remark 2.6} Consider the isothermic surfaces
associated to the cylinder parametrized by
(\ref{cilindroribaucour}). Then from (\ref{eqf}), the functions
$f$ and $g$ are given by
\begin{eqnarray}\label{f}
f=\left\{\begin{array}{ll}
&a_1\cosh(\sqrt{c}\hspace{0,06cm}u_1)+b_1\sinh(\sqrt{c}\hspace{0,06cm}u_1)-\frac{b}{c},\hspace{0,3cm} if \hspace{0,2cm}c>0\\
&a_1\cos(\sqrt{-c}\hspace{0,06cm}u_1)+b_1\sin(\sqrt{-c}\hspace{0,06cm}u_1)-\frac{b}{c},\hspace{0,3cm}
if \hspace{0,2cm}c<0,
\end{array}
\right.
\end{eqnarray}
\begin{eqnarray}\label{g}
g=\left\{\begin{array}{lll}
&a_2\cos(\sqrt{1+c}\hspace{0,06cm}u_2)+b_2\sin(\sqrt{1+c}\hspace{0,06cm}u_2)+\frac{b}{1+c},\hspace{0,3cm} if \hspace{0,2cm}c>-1,c\neq0\\
&\frac{b}{2}u_2^2+a_2u_2+b_2,\hspace{0,3cm} if \hspace{0,2cm}c=-1\\
&a_2\cosh(\sqrt{-1-c}\hspace{0,06cm}u_2)+b_2\sinh(\sqrt{-1-c}\hspace{0,06cm}u_2)+\frac{b}{1+c},\hspace{0,3cm} if \hspace{0,2cm}c<-1\\
\end{array}
\right.
\end{eqnarray}
and from (\ref{algcond}), the constant satisfy the algebric
relation
\begin{eqnarray}
\frac{b^2}{c\big(1+c\big)}=c\big(a_1^2-b_1^2\big)-(1+c)\big(a_2^2+b_2^2\big)\hspace{0,3cm}
&if& \hspace{0,2cm}c>0,\label{a1b1}\\
\frac{b^2}{c\big(1+c\big)}=c\big(a_1^2+b_1^2\big)-(1+c)\big(a_2^2+b_2^2\big)\hspace{0,3cm}
&if& \hspace{0,2cm}-1<c<0,\label{a1}\\
a_1^2+b_1^2+a_2^2+b_2^2-\big(b+b_2\big)^2=0\hspace{0,3cm}
&if& \hspace{0,2cm}c=-1,\label{b2}\\
\frac{b^2}{c\big(1+c\big)}=c\big(a_1^2+b_1^2\big)-(1+c)\big(a_2^2-b_2^2\big)\hspace{0,3cm}
&if& \hspace{0,2cm}c<-1.\label{a2b2}
\end{eqnarray}

Moreover, if $\widetilde{X}$ is not $\frac{-1}{2}-cmc$, then
from Corollary 2.5, using (\ref{a1b1}) and (\ref{a2b2}),
respectively, we have $a_1^2-b_1^2>0$ if $c>0$ and $a_2^2-b_2^2>0$
if $c<-1$.

\vspace{.2in}

Using the Remark 2.5, we get immediately

\vspace{.1in}

\noindent {\bf Corollary 2.7} \textit{ Consider the isothermic
surfaces associated to the cylinder parametrized by
(\ref{cilindroribaucour}). Excluding the  $\frac{-1}{2}-cmc$, we
have that up to rigid motions of $R^3$, the surface
$\widetilde{X}$ is determined by the functions}\\

 \textit{ i) If $c>0$, then
\begin{eqnarray}\label{ff1}
f=\sqrt{A_1}\cosh(\sqrt{c}\hspace{0,06cm}u_1)-\frac{b}{c},\hspace{0,5cm}g=\sqrt{B_1}\sin(\sqrt{1+c}\hspace{0,06cm}u_2)+\frac{b}{1+c},
\end{eqnarray}
where $\displaystyle{\frac{b^2}{c(1+c)}=cA_1-(1+c)B_1}$, with $B_1>0$.}\\

 \textit{ ii) If $-1<c<0$, then
\begin{eqnarray}\label{ff2}
f=\sqrt{A_1}\sin(\sqrt{-c}\hspace{0,06cm}u_1)-\frac{b}{c},\hspace{0,5cm}g=\sqrt{B_1}\sin(\sqrt{1+c}\hspace{0,06cm}u_2)+\frac{b}{1+c},
\end{eqnarray}
where $\displaystyle{\frac{b^2}{c(1+c)}=cA_1-(1+c)B_1}$, with
$A_1>0$ and $B_1>0$.} \\

 \textit{ iii) If $c<-1$, then
\begin{eqnarray}\label{ff3}
f=\sqrt{A_1}\sin(\sqrt{-c}\hspace{0,06cm}u_1)-\frac{b}{c},\hspace{0,5cm}g=\sqrt{B_1}\cosh(\sqrt{-1-c}\hspace{0,06cm}u_2)+\frac{b}{1+c},
\end{eqnarray}
where $\displaystyle{\frac{b^2}{c(1+c)}=cA_1-(1+c)B_1}$, with
$A_1>0$.}\\

 \textit{ iv) If $c=-1$, then
\begin{eqnarray}\label{ff4}
\hspace{0,1cm}f=\sqrt{A_1}\sin(u_1)+b,\hspace{0,5cm}g=\frac{b}{2}u_2^2+a_2u_2+b_2,
\end{eqnarray}
where $A_1+a_2^2+b_2^2-\big(b+b_2\big)^2=0$, with $A_1>0$.
}

\vspace{.1in}

\noindent {\bf Proof:} Consider the isothermic surfaces associated
to the cylinder parametrized by
(\ref{cilindroribaucour}) that is not  $\frac{-1}{2}-cmc$. \\
If $c>0$, then using (\ref{f}), we get
\begin{eqnarray}\nonumber
f=\sqrt{A_1}\cosh(\sqrt{c}\hspace{0,06cm}u_1+A_2)-\frac{b}{c},\hspace{0,5cm}g=\sqrt{B_1}\sin(\sqrt{1+c}\hspace{0,06cm}u_2+B_2)+\frac{b}{1+c},\nonumber
\end{eqnarray}
where $\displaystyle{\frac{b^2}{c(1+c)}}=cA_1-(1+c)B_1$, with
$B_1>0$ and  the constants $A_2$ and $B_2$, without loss of
generality, my be considered to be zero. One can verify that the
surfaces with different values of $A_2$, $B_2$ are congruent by
rigid motions of $R^3$. In fact, using the notation
$\widetilde{X}_{bcA_2B_2}$ for the surface $\widetilde{X}$
with fixed constants $A_2$ and $B_2$, we have
\begin{eqnarray}\nonumber
\widetilde{X}_{bcA_2B_2}=R_{\frac{-B_2}{\sqrt{1+c}}}\widetilde{X}_{bc00}\circ
h+T_{\frac{-A_2}{\sqrt{c}}},\nonumber
\end{eqnarray}
where
$h(u_1,u_2)=\displaystyle{\bigg(u_1+\frac{A_2}{\sqrt{c}},u_2+\frac{B_2}{\sqrt{1+c}}\bigg)}$.

Finally, with analogous argument, we have (\ref{ff2}), (\ref{ff3}) and (\ref{ff4}).

\vspace{.2in}

\noindent {\bf Proposition 2.8}\textit{ Consider the isothermic
surfaces associated to the cylinder parametrized by
(\ref{cilindroribaucour}), that is not  $\frac{-1}{2}-cmc$. Let $M=2b+c(f-g)$ and suppose there is $p_0=(u_1^0,u_2^0)$ such that $M(p_0)=0$.
Then $p_0=\displaystyle{\bigg(0,\frac{(4k+1)\pi}{2\sqrt{1+c}}\bigg)}$, 
$k\in Z$ if $c>0$, $p_0=\displaystyle{\bigg(\frac{(2k-1)\pi}{2\sqrt{-c}},\frac{(2k-1)\pi}{2\sqrt{1+c}}\bigg)}$, $k\in Z$  if $-1<c<0$, $p_0=\displaystyle{\bigg(\frac{(4k+1)\pi}{2\sqrt{-c}},0\bigg)}$, $k\in Z$  if $c<-1$ and $p_0=\displaystyle{\bigg(\frac{(2k-1)\pi}{2},0\bigg)}$, $k\in Z$  if $c=-1$. Moreover, the functions $f$ and $g$ given by (\ref{ff1}) and (\ref{ff2})  become
\begin{eqnarray}\label{fp0}
f=\left\{\begin{array}{lll}
&\frac{1}{c}\bigg(\cosh(\sqrt{c}\hspace{0,06cm}u_1)+1\bigg),\hspace{0,3cm} if \hspace{0,2cm}c>0\\
&\frac{-1}{c}\bigg(\sin(\sqrt{-c}\hspace{0,06cm}u_1)+\epsilon_1\bigg),\hspace{0,3cm}
if \hspace{0,2cm}-1\leq c<0, \hspace{0,2cm}\epsilon_1^2=1,\\
&\frac{-1}{c}\bigg(\sin(\sqrt{-c}\hspace{0,06cm}u_1)+1\bigg),\hspace{0,3cm} if \hspace{0,2cm}c<-1
\end{array}
\right.
\end{eqnarray}
\begin{eqnarray}\label{gp0}
g=\left\{\begin{array}{llll}
&\frac{1}{1+c}\bigg(\sin(\sqrt{1+c}\hspace{0,06cm}u_2)-1\bigg) if \hspace{0,2cm}c>0\\
&\frac{1}{1+c}\bigg(\sin(\sqrt{1+c}\hspace{0,06cm}u_2)+\epsilon_1\bigg) if \hspace{0,2cm}-1<c<0, \hspace{0,2cm}\epsilon_1^2=1,\\
&\frac{\epsilon_1}{2}u_2^2,\hspace{0,3cm}
if \hspace{0,2cm}c=-1, \hspace{0,2cm}\epsilon_1^2=1,\\
&\frac{-1}{1+c}\bigg(\cosh(\sqrt{-1-c}\hspace{0,06cm}u_2)-1\bigg) if \hspace{0,2cm}c<-1
\end{array}
\right.
\end{eqnarray}
}

\vspace{.1in}

\noindent {\bf Proof:} Let $p_0=(u_1^0,u_2^0)$ such that $M(p_0)=0$. Using (\ref{eqf}) and (\ref{eqf1''}), in  $p_0$ we have
\begin{eqnarray}\label{fgp0}
f'(u_1^0)=0, \hspace{0,3cm}g'(u_2^0)=0,\hspace{0,3cm} g(u_2^0)=0\hspace{0,3cm} and \hspace{0,3cm}f(u_1^0)=\frac{-2b}{c}.
\end{eqnarray}
If $c>0$, then using (\ref{ff1}) and (\ref{fgp0}), we obtain
$u_1^0=0$, $\sqrt{A_1}=\displaystyle{\frac{-b}{c}}$, $b<0$.\\
Moreover we have $\cos(\sqrt{1+c}\hspace{0,06cm}u_2^0)=0$ and $\sqrt{B_1}\sin(\sqrt{1+c}\hspace{0,06cm}u_2^0)=\displaystyle{\frac{-b}{1+c}}$. Thus $u_2^0=\displaystyle{\frac{(4k+1)\pi}{2\sqrt{1+c}}}$, $k\in Z$ and $\sqrt{B_1}=\displaystyle{\frac{-b}{1+c}}$.
Therefore $f$ and $g$ given by (\ref{ff1}) can be rewritten by
\begin{eqnarray}\label{fgcomb0}
f=\frac{-b}{c}\bigg(\cosh(\sqrt{c}\hspace{0,06cm}u_1)+1\bigg),\hspace{0,5cm}g=\frac{-b}{1+c}\bigg(\sin(\sqrt{1+c}\hspace{0,06cm}u_2)-1\bigg).
\end{eqnarray}
Without loss of generality, we can consider $b=-1$. In fact, substituting $f$ and $g$ above in (\ref{metrica}), (\ref{aa1}) and (\ref{aa2}), we obtain that the first and second fundamental forms of the
$\widetilde{X}$, do not depend on $b$. Thus, $f$ and $g$, are given by (\ref{fp0}) and and (\ref{gp0}).

If $-1<c<0$, then using (\ref{ff2}) and (\ref{fgp0}), we get $u_1^0=\displaystyle{\frac{(2k-1)\pi}{2\sqrt{-c}}}$, and $u_2^0=\displaystyle{\frac{(2k-1)\pi}{2\sqrt{1+c}}}$, $k\in Z$ .\\
If $b<0$, then we have $\sqrt{A_1}=\displaystyle{\frac{b}{c}}$, and $\sqrt{B_1}=\displaystyle{\frac{-b}{1+c}}$. Thus, $f$ and $g$ given by (\ref{ff2}) can be rewritten by
\begin{eqnarray}\label{fgcomb00}
f=\frac{b}{c}\bigg(\sin(\sqrt{-c}\hspace{0,06cm}u_1)-1\bigg),\hspace{0,5cm}g=\frac{-b}{1+c}\bigg(\sin(\sqrt{1+c}\hspace{0,06cm}u_2)-1\bigg).
\end{eqnarray}
 Substituting $f$ and $g$ above in (\ref{metrica}), (\ref{aa1}) and (\ref{aa2}), we obtain that the first and second fundamental forms of the
$\widetilde{X}$, do not depend on $b$. Therefore, without loss of generality,  we can consider $b=-1$. Thus, $f$ and $g$, are given by (\ref{fp0}) and (\ref{gp0}), with $\epsilon_1=-1$.\\
On the other hand, if $b>0$, then  we have $\sqrt{A_1}=\displaystyle{\frac{-b}{c}}$, and $\sqrt{B_1}=\displaystyle{\frac{b}{1+c}}$. Thus, $f$ and $g$ given by (\ref{ff2}) can be rewritten by
\begin{eqnarray}\label{fgcomb00}
f=\frac{-b}{c}\bigg(\sin(\sqrt{-c}\hspace{0,06cm}u_1)+1\bigg),\hspace{0,5cm}g=\frac{b}{1+c}\bigg(\sin(\sqrt{1+c}\hspace{0,06cm}u_2)+1\bigg).
\end{eqnarray}
As before, without loss of generality, we can consider $b=1$. Thus, $f$ and $g$, are given by (\ref{fp0}) and (\ref{gp0}), with $\epsilon_1=1$.

If $c<-1$, then using (\ref{ff3}) and (\ref{fgp0}), we get $u_2^0=0$, $\sqrt{B_1}=\displaystyle{\frac{-b}{c}}$, $b>0$.\\
Moreover we have $\cos(\sqrt{-c}\hspace{0,06cm}u_1^0)=0$ and $\sqrt{A_1}\sin(\sqrt{-c}\hspace{0,06cm}u_1^0)=\displaystyle{\frac{-b}{c}}$. Thus $u_1^0=\displaystyle{\frac{(4k+1)\pi}{2\sqrt{-c}}}$, $k\in Z$ and $\sqrt{B_1}=\displaystyle{\frac{-b}{c}}$.
Therefore $f$ and $g$ given by (\ref{ff3}) can be rewritten by
\begin{eqnarray}\label{fgcombmenos}
f=\frac{-b}{c}\bigg(\sin(\sqrt{-c}\hspace{0,06cm}u_1)+1\bigg),\hspace{0,5cm}g=\frac{-b}{1+c}\bigg(\cosh(\sqrt{-1-c}\hspace{0,06cm}u_2)-1\bigg).
\end{eqnarray}
As before, without loss of generality, we can consider $b=1$. Thus, $f$ and $g$, are given by (\ref{fp0}) and (\ref{gp0}).

If $c=-1$, similarly to the previous cases, using (\ref{ff4}) and (\ref{fgp0}), we have $p_0=\displaystyle{\bigg(\frac{(2k-1)\pi}{2},\frac{-a_2}{b}\bigg)}$ and
\begin{eqnarray}\label{fgcombmenosum}
f=|b|\big(\sin(u_1)+\epsilon_1\big),\hspace{0,5cm}\epsilon_1^2=1,\hspace{0,5cm}g=\frac{1}{2b}\big(bu_2+a_2\big)^2.
\end{eqnarray}
Without loss of generality, we can consider $a_2=0$. In fact, using (\ref{rotacao}), we have
\begin{eqnarray}\nonumber
\widetilde{X}_{bca_2}=R_{\frac{-a_2}{b}}\widetilde{X}_{bc0}\circ
h,\nonumber
\end{eqnarray}
where
$h(u_1,u_2)=\displaystyle{\bigg(u_1,u_2+\frac{a_2}{b}\bigg)}$ and $\widetilde{X}_{bca_2}$ is the surface $\widetilde{X}$
with fixed constant $a_2$. \\
As before, without loss of generality, we can consider $|b|=1$. Thus, $f$ and $g$, are given by (\ref{fp0}) and (\ref{gp0}).

\vspace{.2in}

\noindent {\bf Proposition 2.9} \textit{Any isothermic surfaces
associated to the cylinder $\widetilde{X}$, given by
Theorem 2.1 is complete.}

\vspace{.1in}

\noindent {\bf Proof:}
For any divergent curve $\gamma:[0,1)\rightarrow V-\{p_0\} $ such that $\lim_{t\rightarrow 1}\gamma(t)=p_0$ where $p_0$ are given by Proposition 2.8, then from Proposition 2.5 the length of $\widetilde{X}(\gamma)$
is infinite.\\
For divergent curves
$\gamma(t)=(u_1(t),u_2(t))$, such that
$\displaystyle{\lim_{t\rightarrow
\infty}\big(u_1^2+u_2^2\big)=\infty}$, we have
$l(\widetilde{X}\circ\gamma)=\infty$.

In fact, if $c>0$, then the functions $f$ and $g$ are given by
(\ref{ff1}), and from (\ref{metrica}) the coefficients of the
first fundamental form is
$\psi=\displaystyle{\frac{|c\big(f+g\big)|}{|2b+c\big(f-g\big)|}}$.
Therefore,  $\displaystyle{\lim_{|u_1|\rightarrow \infty}\psi=1}$
uniformly in $u_2$. Hence, there exist $k>0$ such that $|\psi(u_1,u_2)|>\frac{1}{2}$ for $(u_1,u_2)\in
R^2$ with $|u_1|>k$. Let
\begin{eqnarray}
m=min\bigg\{|\psi(u_1,u_2)|;\hspace{0,1cm}(u_1,u_2)\in
\hspace{0,05cm}[-k,k]\times
\bigg[0,\frac{2\pi}{\sqrt{1+c}}\bigg]\bigg\}.
\end{eqnarray}
Note that, $g(u_2)=g(u_2+\frac{2\pi}{\sqrt{1+c}})$, therefore
$|\psi(u_1,u_2)|>m$ in $[-k,k]\times R$. Consider\\
$m_0=min\{m,\frac{1}{2}\}$, then $|\psi(u_1,u_2)|>m_0$ in $R^2$.
Thus $l(\widetilde{X}\circ\gamma)=\infty$. The case $c\leq-1$
is analogous.

  Finally if $-1<c<0$, then the functions $f$ and $g$ are given by
(\ref{ff2}), and from (\ref{metrica}) the coefficients of the
first fundamental form is
$\psi=\displaystyle{\frac{|c\big(f+g\big)|}{|2b+c\big(f-g\big)|}}$.
In this case, let
\begin{eqnarray}
m_0=min\bigg\{|\psi(u_1,u_2)|;\hspace{0,1cm}(u_1,u_2)\in
\hspace{0,05cm}\bigg[0,\frac{2\pi}{\sqrt{-c}}\bigg]\times
\bigg[0,\frac{2\pi}{\sqrt{1+c}}\bigg]\bigg\}.
\end{eqnarray}
Note that, $g(u_2)=g(u_2+\frac{2\pi}{\sqrt{-1-c}})$ and  $f(u_1)=f(u_1+\frac{2\pi}{\sqrt{-c}})$, therefore
$|\psi(u_1,u_2)|>m_0$ in $R^2$.
Thus $l(\widetilde{X}\circ\gamma)=\infty$ and we conclude
that $\widetilde{X}$ is a complete surface.

\vspace{.15in}

\noindent {\bf Remark 2.10} As obtained in \cite{8}, the $n$ points of maximum ( respectively minimum ) for the Gaussian curvature of the family of
isothermic surfaces given by (\ref{cilindroribaucour}), generate 1-periodic n-bubble surfaces outside of the cylinder. In addition, we can also have 1-periodic n-bubble surfaces inside of the cylinder, since $2b+c(f-g)\neq 0$. We see this with some examples. 

\vspace{.1in}

i) For $c>0$, if $\sqrt{1+c}=\displaystyle{\frac{n}{m}}$, then
generate 1-periodic n-bubble surfaces with two ends of geometric
index m . And more, if $b>0$ then the n-bubble surfaces inside ( see Figure 1 ),
if $b<0$ them the n-bubble surfaces outside ( see Figure 2 ).


\begin{figure}[h]

\centering
 \includegraphics[scale=0.4]{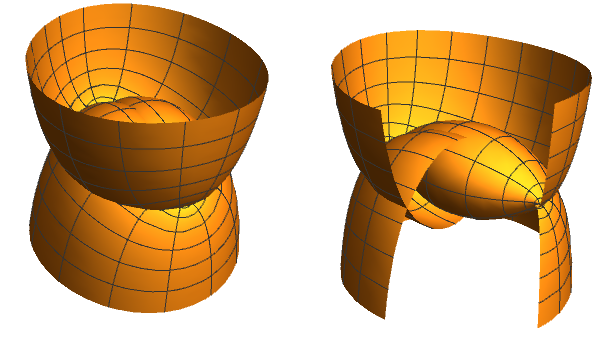}
\label{figura1}
\caption{In the Figure 1 above, we have $c=3$, $b=4\sqrt{6}$ and $\sqrt{1+c}=\displaystyle{\frac{2}{1}}$. Thus  we have 2-bubble surfaces inside.}
\end{figure}

\vspace{.9in}
\vspace{.9in}

\vspace{.2in}

\begin{figure}[h]

\centering
\includegraphics[scale=0.5]{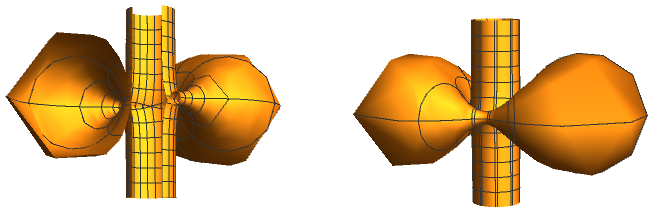}
\caption{In the Figure 2 above we have $c=3$, $b=-4\sqrt{6}$ and $\sqrt{1+c}=\displaystyle{\frac{2}{1}}$. Thus we have 2-bubble surfaces ouside.}

\end{figure}

\vspace{.2in}

ii) For $-1<c<0$, if  $\sqrt{-c}$ and $\sqrt{1+c}$ are a rational number, then we have  1-periodic bubbles surfaces inside  1-periodic bubbles surfaces vertically ( see Figure 3 ). If $\sqrt{-c}$ and (or) $\sqrt{1+c}$ are not rational number, then generate bubble surfaces inside bubbles surfaces, but not 1-periodic.

\vspace{.2in}

\begin{figure}[h]

\centering
\includegraphics[scale=0.5]{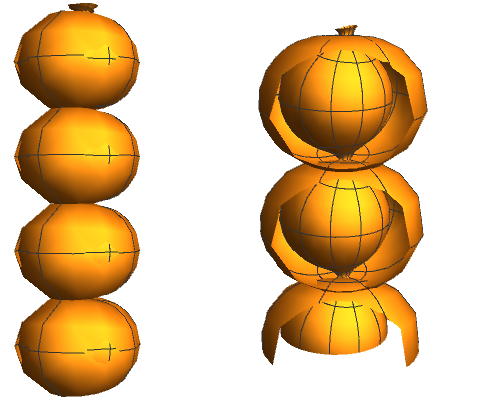}
\caption{In the Figure 3 above we have $c=\frac{-16}{25}$,  $b=\frac{12\sqrt{73}}{125}$ and $\sqrt{-c}=\displaystyle{\frac{4}{5}}$ $\sqrt{1+c}=\displaystyle{\frac{3}{5}}$.}

\end{figure}

\vspace{.2in}

\vspace{.2in}

iii) For $c<-1$, if $b>0$, then we have 1-periodic bubbles surfaces outside ( see Figure 4 ). If $b<0$, then we have 1-periodic bubbles surfaces inside ( see Figure 5 ).

\begin{figure}[h]

\centering
\includegraphics[scale=0.4]{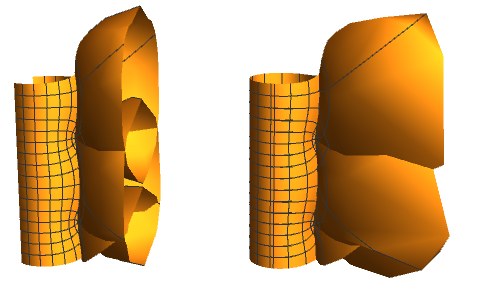}
\caption{In the Figure 4 above we have $c=-5$,  $b=\frac{4\sqrt{5}}{3}>0$}

\end{figure}

\begin{figure}[h]

\centering
\includegraphics[scale=0.4]{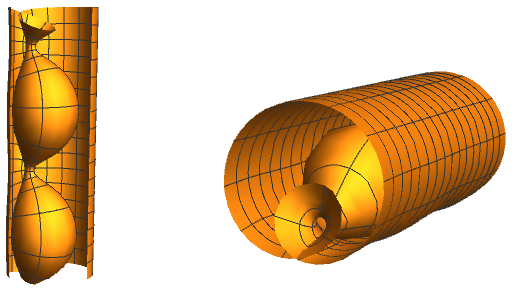}
\caption{In the Figure 5 above we have $c=-5$,  $b=\frac{-4\sqrt{5}}{3}<0$}

\end{figure}

iv)  For $c=-1$, in this case we have 1-periodic bubble surfaces outside vertically ( see Figure 6 )

\begin{figure}[h]

\centering
\includegraphics[scale=0.4]{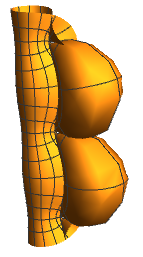}
\caption{In Figure 6 the left above, we have $b=2$. }

\end{figure}

\vspace{.2in}

\noindent {\bf Remark 2.11} As obtained in \cite{9}, the family of
isothermic surfaces given by (\ref{cilindroribaucour}), generate planar ends. This occurs whenever there are $p_0\in R^2$ such that $\big(2b+c(f-g)\big)(p_0)=0$. In this case, $p_0$ are given by Proposition 2.8 and the functions $f$ and $g$ are given by (\ref{fp0}) and (\ref{gp0}). Besides that, the $\widetilde{X}$ depends only on a parameter.

\vspace{.1in}

\begin{figure}[h]

\centering
\includegraphics[scale=0.6]{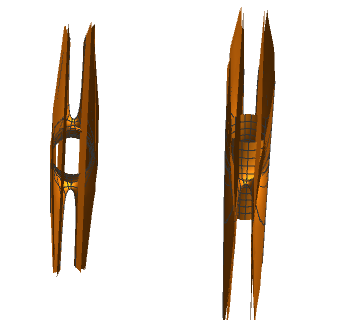}
\caption{In Figure 7 above, we have $c=3$, in $p_0=\bigg(0,\frac{(4k+1)\pi}{4}\bigg)$, $k\in Z$ we have planar ends. }

\end{figure}

\vspace{.1in}

\begin{figure}[h]

\centering
\includegraphics[scale=0.6]{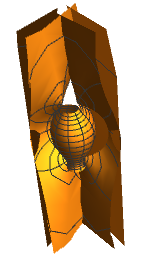}
\caption{In Figure 8 above, we have $c=\frac{-16}{25}$, in $p_0=\bigg(\frac{5(2k-1)\pi}{8},\frac{5(2k-1)\pi}{6}\bigg)$ we have planar ends. }

\end{figure}

\vspace{.1in}

\begin{figure}[h]

\centering
\includegraphics[scale=0.6]{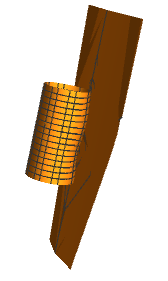}
\caption{In Figure 9 above, we have  $c=-5$, in $p_0=\bigg(\frac{(4k+1)\pi}{2\sqrt{5}},0\bigg)$ we have planar ends . }

\end{figure}

\section*{ \large \bf 3. Solution of the Calapso Equation.}

\vspace{.2in}

In \cite{3}, we have that for each isothermic surface, we have a solution of the Calapso equation. In the section, we obtain solution of the Calapso Equation.

The Calapso equation defined in [Calapso] given by
\begin{eqnarray}\label{calapso}
\bigg(\frac{\omega,_{u_1u_2}}{w}\bigg),_{u_1u_1}+\bigg(\frac{\omega,_{u_1u_2}}{w}\bigg),_{u_2u_2}+\big(\omega^2\big),_{u_1u_2}=0,
\end{eqnarray}
describes isoghermic surfaces in $R^3$.

\vspace{.1in}

 \noindent {\bf Remark 3.1} \cite{13} Let $X(u_1,u_2)$ be an isothermic surface with the first fundamental form givem by
 \begin{eqnarray}\nonumber
I=e^{2\varphi}\big(du_1^2+du_2^2\big).\nonumber
\end{eqnarray}
 Then the functions $\omega=\sqrt{2}e^{\varphi}H$ and $\Omega=\sqrt{2}e^{\varphi}H'$ are solutions of the Calapso equation, where $H$ is the mean curvature of $X$ and $H'$ is the skew curvature of $M$.

 \vspace{.2in}

\noindent {\bf Proposition 3.2} \textit{Consider the isothermic
surfaces associated to the cylinder parametrized by
(\ref{cilindroribaucour}). Then the mean and skew curvature of the
$\widetilde{X}$ are given by
\begin{eqnarray}
\widetilde{H}&=&\frac{-1}{2}-\frac{b}{c\big(f+g\big)},\label{meancurv}\\
\widetilde{H}'&=&\frac{M-2b}{M\psi^2},\label{skewcurv}
\end{eqnarray}
where the functions $f$ and $g$ are given by (\ref{ff1})-(\ref{ff4}), $M=2b+c(f-g)$ and $\psi$ are the coefficients of the first fundamental form of the  $\widetilde{X}$, given by (\ref{metrica}).
}

\vspace{.1in}

\noindent {\bf Proof:} In fact, the (\ref{meancurv}) is given by Proposition 2.3. For provide (\ref{skewcurv}), just replace (\ref{aa1}) and (\ref{aa2}) in $\widetilde{H}'=\frac{\widetilde{\lambda}_1-\widetilde{\lambda}_2}{2}$ and we conclude the prove.\\

 \vspace{.2in}
 
 Using Remark 3.1 and Proposition 3.2, we get immediately
 
  \vspace{.1in}
 
  \noindent {\bf Corollary 3.3} \textit{Consider the isothermic
surfaces associated to the cylinder parametrized by
(\ref{cilindroribaucour}), whose first fundamental form is given by 
\begin{eqnarray}\nonumber
I=\bigg(\frac{c\big(f+g\big)}{M}\bigg)^2\big(du_1^2+du_2^2\big),\nonumber
\end{eqnarray}
where the functions $f$ and $g$ are given by (\ref{ff1})-(\ref{ff4}) and $M=2b+c(f-g)$. \\
 Then the functions $\omega=\displaystyle{\frac{\epsilon\sqrt{2}\big(M+2cg\big)}{2M}}$ and $\Omega=\displaystyle{\frac{\epsilon\sqrt{2}(f-g)}{f+g}}$, with $\epsilon=1$ if $c>0$ and  $\epsilon=-1$ if $c<0$, are solutions of the Calapso equation.}

 \vspace{.2in}
 
  \noindent {\bf Example 3.4} \textit{Consider the isothermic
surfaces associated to the cylinder given by Figure 1.\\
In this case, we have $b=4\sqrt{6}$, $c=3$, $f(u_1)=2\cosh(\sqrt{3}u_1)-\frac{4\sqrt{6}}{3}$ and $g(u_2)=\sin(2u_2)+\sqrt{6}$. \\
Using the Corollary 3.3, we have
\begin{eqnarray}
\omega=\frac{\sqrt{2}\big(7\sqrt{6}+6\cosh(\sqrt{3}u_1)+3\sin(2u_2)\big)}{2\sqrt{6}+12\cosh(\sqrt{3}u_1)-6\sin(2u_2)}\nonumber\\
\Omega=\frac{\sqrt{2}\big(-7\sqrt{6}+6\cosh(\sqrt{3}u_1)-3\sin(2u_2)\big)}{-\sqrt{6}+6\cosh(\sqrt{3}u_1)+3\sin(2u_2)}\nonumber
\end{eqnarray}}

\begin{figure}[h]

\centering
\includegraphics[scale=0.5]{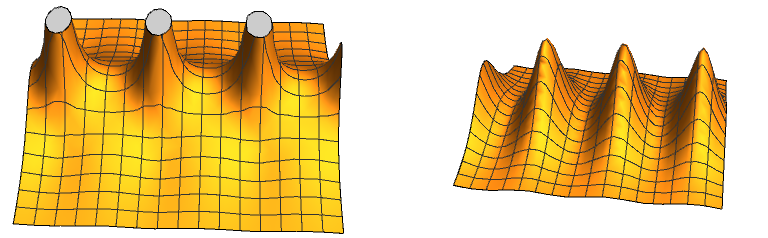}
\caption{In the figure above we have the graphics of the solutions of the Calapso equation.}

\end{figure}

\vspace{.2in}
 
  \noindent {\bf Example 3.5} \textit{Consider the isothermic
surfaces associated to the cylinder given by Figure 2.\\
In this case, we have $b=-4\sqrt{6}$, $c=3$, $f(u_1)=2\cosh(\sqrt{3}u_1)+\frac{4\sqrt{6}}{3}$ and $g(u_2)=\sin(2u_2)-\sqrt{6}$. \\
Using the Corollary 3.3, we have
\begin{eqnarray}
\omega=\frac{\sqrt{2}\big(-7\sqrt{6}+6\cosh(\sqrt{3}u_1)+3\sin(2u_2)\big)}{-2\sqrt{6}+12\cosh(\sqrt{3}u_1)-6\sin(2u_2)}\nonumber\\
\Omega=\frac{\sqrt{2}\big(7\sqrt{6}+6\cosh(\sqrt{3}u_1)-3\sin(2u_2)\big)}{\sqrt{6}+6\cosh(\sqrt{3}u_1)+3\sin(2u_2)}\nonumber
\end{eqnarray}}
The graphics for these solutions of the Calapso equation are similar to the graphics in Figure 6.
 
  \vspace{.2in}

  \noindent {\bf Example 3.6} \textit{Consider the isothermic
surfaces associated to the cylinder given by Figure 3.\\
In this case, we have $c=\frac{-16}{25}$,  $b=\frac{12\sqrt{73}}{125}$, $f(u_1)=2\sin(\frac{4}{5}u_1)+\frac{3\sqrt{73}}{20}$ and $g(u_2)=\sin(\frac{3}{5}u_2)+\frac{4\sqrt{73}}{15}$. \\
Using the Corollary 3.3, we have
\begin{eqnarray}
\omega=\frac{\sqrt{2}\big(7\sqrt{73}+120\sin(\frac{4}{5}u_1)+60\sin(\frac{3}{5}u_2)\big)}{50\sqrt{73}-240\sin(\frac{4}{5}u_1)+120\sin(\frac{3}{5}u_2)}\nonumber\\
\Omega=\frac{\sqrt{2}\big(7\sqrt{73}-120\sin(\frac{4}{5}u_1)+60\sin(\frac{3}{5}u_2)\big)}{25\sqrt{73}+120\sin(\frac{4}{5}u_1)+60\sin(\frac{3}{5}u_2)}\nonumber
\end{eqnarray}}

\begin{figure}[h]

\centering
\includegraphics[scale=0.5]{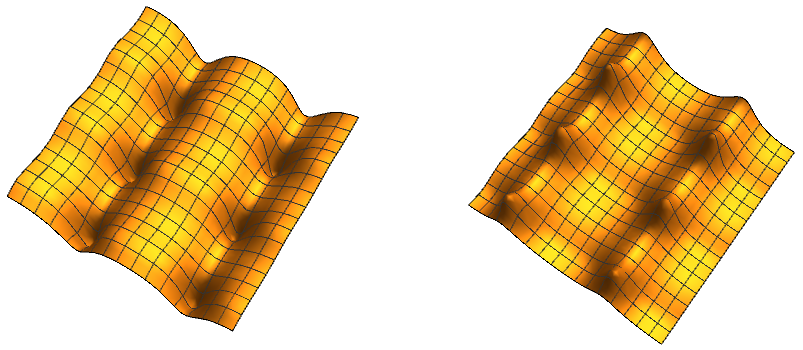}
\caption{In the figure above we have the graphics of the Example 3.6 solutions of the Calapso equation .}

\end{figure}

 \vspace{.2in}

  \noindent {\bf Example 3.7} \textit{Consider the isothermic
surfaces associated to the cylinder given by Figure 4.\\
In this case, we have $c=-5$,  $b=\frac{4\sqrt{5}}{3}$, $f(u_1)=\frac{\sin(\sqrt{5}u_1)}{3}+\frac{4\sqrt{5}}{15}$ and $g(u_2)=\frac{\cosh(2u_2)}{2}-\frac{\sqrt{5}}{3}$. \\
Using the Corollary 3.3, we have
\begin{eqnarray}
\omega=\frac{-\sqrt{2}\big(-18\sqrt{5}+10\sin(\sqrt{5}u_1)+15\cosh(2u_2)\big)}{4\sqrt{5}+20\sin(\sqrt{5}u_1)-30\cosh(2u_2)}\nonumber\\
\Omega=\frac{-\sqrt{2}\big(18\sqrt{5}+10\sin(\sqrt{5}u_1)-15\cosh(2u_2)\big)}{-2\sqrt{5}+10\sin(\sqrt{5}u_1)+15\cosh(2u_2)}\nonumber
\end{eqnarray}}

\begin{figure}[h]

\centering
\includegraphics[scale=0.4]{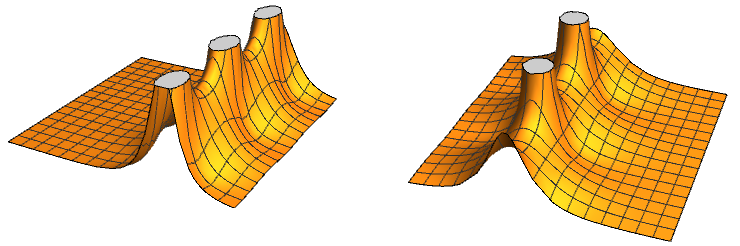}
\caption{In the figure above we have the graphics of the Example 3.7 solutions of the Calapso equation .}

\end{figure}

\vspace{.2in}

  \noindent {\bf Example 3.8} \textit{Consider the isothermic
surfaces associated to the cylinder given by Figure 6.\\
In this case, we have $c=-1$,  $b=2$, $f(u_1)=\sin(u_1)+2$ and $g(u_2)=u_2^2-\frac{\sqrt{3}}{4}$. \\
Using the Corollary 3.3, we have
\begin{eqnarray}
\omega=\frac{-\sqrt{2}\big(-11+4u_2^2+4\sin(u_1)\big)}{-10-8u_2^2+8\sin(u_1)}\nonumber\\
\Omega=\frac{-\sqrt{2}\big(11-4u_2^2+4\sin(u_1)\big)}{5+4u_2^2+4\sin(u_1)}\nonumber
\end{eqnarray}}

\begin{figure}[h]

\centering
\includegraphics[scale=0.4]{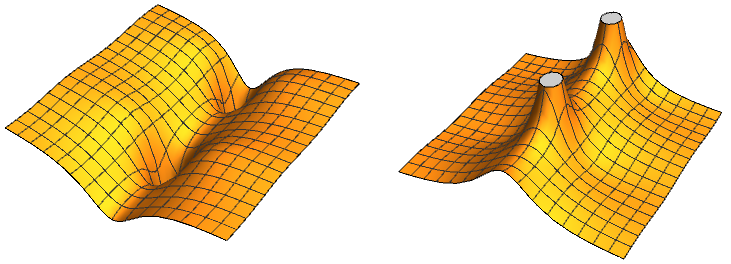}
\caption{In the figure above we have the graphics of the Example 3.8 solutions of the Calapso equation .}

\end{figure}

\vspace{.2in}

  \noindent {\bf Example 3.9} \textit{Consider the isothermic
surfaces associated to the cylinder given by Figure 7.\\
In this case, we have $c=3$, $f(u_1)=\frac{\cosh(u_1)+1}{3}$ and $g(u_2)=\frac{\sin(2u_2)-1}{4}$. \\
Using the Corollary 3.3, we have
\begin{eqnarray}
\omega=\frac{\sqrt{2}\big(7+4\cosh(u_1)-3\sin(2u_2)\big)}{-2+8\cosh(u_1)-6\sin(2u_2)}\nonumber\\
\Omega=\frac{\sqrt{2}\big(-7-4\cosh(u_1)+3\sin(2u_2)\big)}{1+4\cosh(u_1)+3\sin(2u_2)}\nonumber
\end{eqnarray}}

\begin{figure}[h]

\centering
\includegraphics[scale=0.3]{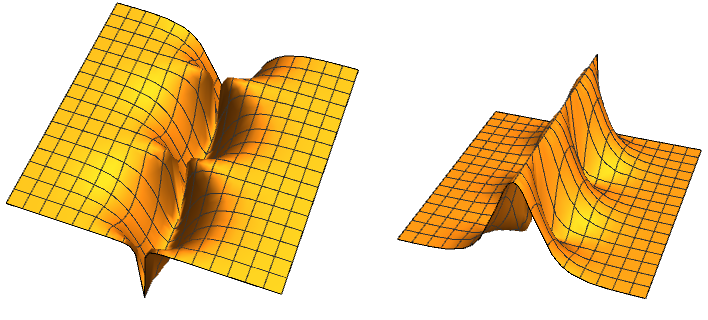}
\caption{In the figure above we have the graphic of the the function $\omega$.}

\end{figure}

\vspace{.2in}

 \noindent {\bf Example 3.10} \textit{Consider the isothermic
surfaces associated to the cylinder given by Figure 8.\\
In this case, we have $c=\frac{-16}{25}$, $f(u_1)=\frac{25\sin(\frac{4}{5}u_1)+25}{16}$ and $g(u_2)=\frac{25\sin(\frac{3}{5}u_2)+25}{16}$. \\
Using the Corollary 3.3, we have the solution  of the Calapso Equation for this case. The graphic are similar to the graphs in the previous examples.}

\vspace{.2in}

 \noindent {\bf Example 3.11} \textit{Consider the isothermic
surfaces associated to the cylinder given by Figure 9.\\
In this case, we have $c=-5$, $f(u_1)=\frac{\sin(\sqrt{5}u_1)+1}{5}$ and $g(u_2)=\frac{\cosh(2u_2)-1}{4}$. \\
Using the Corollary 3.3, we have the solution  of the Calapso Equation for this case. The graphic are similar to the graphs in the previous examples.}
 
 \vspace{.2in}
 
 \noindent {\bf Remark 3.12}
 For each isothermic surface obtained in this work, we can apply isometries, dilations, inversions, obtaining new isothermic surfaces. So we can get new solutions of the Calapso equation.
 {}

 }
\end{document}